\documentclass[12pt]{amsart}
\usepackage{amssymb}
\usepackage{enumerate}
\usepackage{amsthm}
\usepackage{tikz-cd}
\usepackage{amsmath}
\usepackage[mathscr]{euscript}
\DeclareMathOperator{\rank}{{rank}}
\DeclareMathOperator{\Eff}{{Eff}}
\DeclareMathOperator{\Nef}{{Nef}}
\makeatletter
\@namedef{subjclassname@2010}{%
  \textup{2010} Mathematics Subject Classification}
\makeatother
\newtheorem{thm}{Theorem}[section]

\begin{document}

\baselineskip=17pt

\author{Rupam Karmakar}
\address{Institute of Mathematical Sciences\\ CIT Campus, Taramani, Chennai 600113, India and Homi Bhabha National Institute, Training School Complex, Anushakti Nagar, Mumbai 400094, India}
\email[Rupam Karmakar]{rupamk@imsc.res.in}
\begin{abstract}
Let $X = \mathbb{P}(E_1) \times_C \mathbb{P}(E_2)$ where $C$ is a smooth curve and $E_1$, $E_2$ are vector bundles over $C$.In this paper we compute the pseudo effective cones of higher codimension cycles on $X$.
\end{abstract}

\title{Effective cones of cycles on products of projective bundles over curves}
\maketitle
\section{Introduction}
The  cones of divisors and curves on projective varieties have been extensively studied over the years and by now are quite well understood. However, more recently the theory of cones of cycles of higher dimension has been the subject of increasing interest(see \cite{F}, \cite{DELV}, \cite{DJV}, \cite{CC} etc). Lately, there has been significant progress in the theoretical understanding of such cycles, due to \cite{FL1}, \cite{FL2} and others. But the the number of examples where the cone of effective cycles have been explicitly computed is relatively small till date \cite{F}, \cite{CLO}, \cite{PP} etc.

Let $E_1$ and $E_2$ be two vector bundles over a smooth curve $C$ and consider the fibre product  $X = \mathbb{P}(E_1) \times_C \mathbb{P}(E_2)$. Motivated by the results in \cite{F}, in this paper, we compute the cones of effective cycles on $X$ in the following cases.

Case I: When both $E_1$ and $E_2$ are semistable vector bundles of rank $r_1$ and $r_2$ respectively, the cone of effective codimension k-cycles are described in theorem 3.2.

Case II: When Neither $E_1$ nor $E_2$ is semistable, the cone of low dimension effective cycles are computed in theorem 3.3 and the remaining cases in therem 3.5.
\section{Preliminaries}
 Let $X$ be a smooth projective varity of dimension $n$. $N_k(X)$ is the real vector space of k-cycles on $X$  modulo numerical equivalence. For each $k$, $N_k(X)$ is a real vector space of finite dimension. Since $X$ is smooth, we can identify $N_k(X)$ with the abstract dual $N^{n - k}(X) := N_{n - k}(X)^\vee$ via the intersection pairing $N_k(X) \times N_{n - k}(X) \longrightarrow \mathbb{R}$.
 
  For any k-dimensional subvariety $Y$ of $X$, let $[Y]$ be its class in $N_k(X)$. A class $\alpha \in N_k(X)$ is said to be effective if there exists subvarities $Y_1, Y_2, ... , Y_m$ and non-negetive real numbers $n_1, n
 _2, ..., n_m$ such that $\alpha$ can be written as $ \alpha = \sum n_i Y_i$. The \textit{pseudo-effective cone} \,  $\overline{\Eff}_k(X) \subset N_k(X)$ is the closure of the cone generated by classes of effective cycles. It is full-dimensional and does not contain any nonzero linear subspaces. The pseudo-effective dual classes form a closed cone in $N^k(X)$ which we denote by $\overline{\Eff}^k(X)$.
 
 For smooth varities $Y$ and $Z$, a map $f: N^k(Y) \longrightarrow N^K(Z)$ is called pseudo-effective if $f(\overline{\Eff}^k(Y)) \subset \overline{\Eff}^k(Z)$.
 
 The \textit{nef cone} $\Nef^k(X) \subset N^k(X)$ is the dual of $\overline{\Eff}^k(X) \subset N^k(X)$ via the pairing $N^k(X) \times N_k(X) \longrightarrow \mathbb{R}$, i.e,
 \begin{align*}
 \Nef^k(X) := \Big\{ \alpha \in N_k(X) | \alpha \cdot \beta \geq 0 \forall \beta \in \overline{\Eff}_k(X) \Big\}
 \end{align*}
 \section{Cone of effective cycles}
Let $E_1$ and $E_2$ be two vector bundles over a smooth curve $C$ of rank $r_1$, $r_2$ and degrees $d_1$, $d_2$ respectively. Let $\mathbb{P}(E_1) = \bf Proj $ $(\oplus_{d \geq 0}Sym^d(E_1))$ and $\mathbb{P}(E_2) = \bf Proj $ $(\oplus_{d \geq 0}Sym^d(E_2))$ be the 
associated projective  bundle together with the projection morphisms $\pi_1 : \mathbb{P}(E_1) \longrightarrow C$ and  $\pi_2 : \mathbb{P}(E_2) \longrightarrow C$ respectively. Let $X = \mathbb{P}(E_1) \times_C \mathbb{P}(E_2)$ be the fibre product over
$C$. Consider the following commutative diagram:
\begin{center}
\begin{equation}
\begin{tikzcd}
 X = \mathbb{P}(E_1) \times_C \mathbb{P}(E_2) \arrow[r, "p_2"] \arrow[d, "p_1"]
& \mathbb{P}(E_2)\arrow[d,"\pi_2"]\\
\mathbb{P}(E_1) \arrow[r, "\pi_1" ]
& C
\end{tikzcd}
\end{equation}
\end{center}

Let $f_1,f_2$ and $F$ denote the numerical equivalence classes of the fibres of the maps $\pi_1,\pi_2$ and $\pi_1 \circ p_1 = \pi_2 \circ p_2$ respectively. Note that, $X \cong \mathbb{P}(\pi_1^*(E_2)) \cong \mathbb{P}(\pi_2^*(E_1))$.
We first fix the following notations for the numerical equivalence classes,
\vspace{2mm}
 \begin{center}
$\eta_1 = [\mathcal{O}_{\mathbb{P}(E_1)}(1)] \in N^1(\mathbb{P}(E_1))$  \hspace{3mm}, \hspace{3mm} $\eta_2 = [\mathcal{O}_{\mathbb{P}(E_2)}(1)] \in N^1(\mathbb{P}(E_2)),$
\end{center}
\begin{center}
$\xi_1 = [\mathcal{O}_{\mathbb{P}(\pi_1^*(E_2))}(1)]$  \hspace{2mm}, \hspace{2mm} $\xi_2 = [\mathcal{O}_{\mathbb{P}(\pi_2^*(E_1))}(1)]$\hspace{2mm} ,\hspace{2mm} $\zeta_1 = p_1^*(\eta_1)$\hspace{2mm} , \hspace{2mm} $\zeta_2 = p_2^*(\eta_2) $
 
 \vspace{2mm}
 \end{center}
 \begin{center}
 $ \zeta_1 = \xi_2$,\, $\zeta_2 = \xi_1$ \hspace{2mm}, \hspace{2mm} $F= p_1^\ast(f_1) = p_2^\ast(f_2)$
 \end{center}
 \bigskip
 We here summarise some results that has been discussed in \cite{KMR} ( See section 3 in \cite{KMR} for more details) :
 \begin{center}
  $N^1(X)_\mathbb{R} = \mathbb{R}\zeta_1 \oplus \mathbb{R}\zeta_2 \oplus \mathbb{R}F,$
 
$\zeta_1^{r_1}\cdot F = 0\hspace{1.5mm},\hspace{1.5mm}  \zeta_1^{r_1 + 1} = 0 \hspace{1.5mm}, \hspace{1.5mm} \zeta_2^{r_2}\cdot F = 0 \hspace{1.5mm}, \hspace{1.5mm} \zeta_2^{r_2 + 1} = 0 \hspace{1.5mm}, \hspace{1.5mm} F^2 = 0$ ,

$\zeta_1^{r_1} = (\deg(E_1))F\cdot\zeta_1^{r_1-1}\hspace{1.5mm},  \hspace{1.5mm}   \zeta_2^{r_2} = (\deg(E_2))F\cdot\zeta_2^{r_2-1}\hspace{3.5mm}, \hspace{3.5mm}$

$\zeta_1^{r_1}\cdot\zeta_2^{r_2-1} = \deg(E_1)\hspace{3.5mm}, \hspace{3.5mm} \zeta_2^{r_2}\cdot\zeta_1^{r_1-1} = \deg(E_2)$.

 \end{center}

\smallskip
Also, The dual basis of $N_1(X)_\mathbb{R}$ is given by $\{\delta_1, \delta_2, \delta_3\}$, where
 \begin{center}
  
$\delta_1 = F\cdot\zeta_1^{r_1-2}\cdot\zeta_2^{r_2-1}, $

$\delta_2 = F\cdot\zeta_1^{r_1-1}\cdot\zeta_2^{r_2-2},$

 $\delta_3 = \zeta_1^{r_1-1}\cdot\zeta_2^{r_2-1} - \deg(E_1)F\cdot\zeta_1^{r_1-2}\cdot\zeta_2^{r_2-1} - \deg(E_2)F\cdot\zeta_1^{r_1-1}\cdot\zeta_2^{r_2-2}.$
 \end{center}
 
 \begin{thm}
 Let $r_1 = \rank(E_1)$ and $r_2 = \rank(E_2)$ and without loss of generality assume that $ r_1 \leq r_2$. Then the bases of $N^k(X)$ are given by
 
 $$
 N^k(X) = 
  \begin{cases}
 \Big( \{ \zeta_1^i \cdot \zeta_2^{k - i}\}_{i = 0}^k, \{ F \cdot \zeta_1^j \cdot \zeta_2^{k - j - 1} \}_{j = 0}^ {k - 1} \Big) &  if \quad k < r_1\\ \\

 \Big( \{ \zeta_1^i \cdot \zeta_2^{k - i} \}_{i = 0} ^{r_1 - 1}, \{ F \cdot \zeta_1^j \cdot \zeta_2^{k - j - 1} \}_{j = 0}^ {r_1 - 1} \Big) &  if \quad r_1 \leq k < r_2 \\ \\
 
 \Big( \{ \zeta_1^i \cdot \zeta_2^{k - i} \}_{i = t+1} ^{r_1 - 1} , \{ F \cdot \zeta_1^j \cdot \zeta_2^{k - j - 1} \}_{j = t}^{r_1 - 1}  \Big) & if \quad k = r_2 + t \quad where \quad t \in \{0, 1, 2, ..., r_1 - 2 \}.
 
   \end{cases}
 $$
 \end{thm}
 \vspace{4mm}
 \begin{proof}
 To begin with consider the case where $ k < r_1$. We know that $X \cong \mathbb{P}(\pi_2^* E_1)$ and the natural morphism $ \mathbb{P}(\pi_2^*E_1) \longrightarrow \mathbb{P}(E_2)$ can be identified with $p_2$. With the above identifications in place the chow group of $X$ has the following isomorphism [see theorem 3.3 , page 64 \cite{Ful}]
 \begin{align}
 A(X) \cong \bigoplus_{i = 0}^ {r_1 - 1} \zeta_1^i A(\mathbb{P}(E_2))
 \end{align}
\quad Choose $i_1, i_2$ suct that $ 0\leq i_1 < i_2 \leq k$. Consider the $k$- cycle  $\alpha := F \cdot \zeta_1^{r_1 - i_1 -1} \cdot \zeta_2^{r_2 + i_1 -k - 1}$.
 
 Then $\zeta_1^{i_1} \cdot \zeta_2^{k - i_1} \cdot \alpha = 1$ but $ \zeta_1^{i_2} \cdot \zeta_2^{k - i_2} \cdot \alpha = 0$. So, $\{ \zeta_1^{i_1} \cdot \zeta_2^{k - i_1} \}$ and $\{\zeta_1^{i_2} \cdot \zeta_2^{k - i_2} \}$ can not be numerically equivalent.
 
 Similarly, take $j_1, j_2$ such that $ 0 \leq j_1 < j_2 \leq k$ and consider the $k$-cycle\\ $\beta := \zeta_1^{r_1 - j_1 - 1} \cdot \zeta_2^{r_2 + j_1 - k}$.
 
 Then as before it happens that $F \cdot \zeta_1^{j_1} \cdot \zeta_2^{k - j_1 - 1} \cdot \beta = 1$ but $F\cdot \zeta_1^{j_2} \cdot \zeta_2^{k - j_2 - 1} \cdot \beta = 0$. So $\{ F \cdot \zeta_1^{j_1} \cdot \zeta_2^{k - j_1 - 1}\}$ and $\{ F \cdot \zeta_1^{j_2} \cdot \zeta_2^{k - j_2 - 1} \}$ can not be numerically equivalent.
 
 For the remaining case lets assume $0 \leq i \leq j \leq k$ and consider the k-cycle $\gamma := F \cdot \zeta_1^{r_1 - i -1} \cdot \zeta_2^{r_2 + i - 1 - k}$.
 
 Then $\{ \zeta_1^{i} \cdot \zeta_2^{k - i} \} \cdot \gamma = 1$ and $ \{F \cdot \zeta_1^j \cdot \zeta_2^{k - j - 1} \} \cdot \gamma = 0$. So, they can not be numerically equivalent.
 From these observations and $(2)$ we obtain a basis of $N^k(X)$ which is given by
 \begin{align*}
 N^k(X) = \Big( \{ \zeta_1^i \cdot \zeta_2^{k - i} \}_{i = 0}^ k, \{ F \cdot \zeta_1^j \cdot \zeta_2^{k - j - 1} \}_{j = 0}^{k - 1} \Big)
 \end{align*}
 
 For the case $ r_1 \leq k < r_2$  observe that $ \zeta_1^{r_1 + 1} = 0$, $ F \cdot \zeta_1^ {r_1} = 0$ and $ \zeta_1^ {r_2} = deg(E_1)F \cdot \zeta_1^{ r_1 - 1}$.
 
 When $k \geq r_2$  we write as $k = r_2 + t$ where $t$ ranges from  $ 0$ to $r_1 - 1$. In that case the observations like $ \zeta_2^{r_2 + 1} = 0$, $F\cdot\zeta_2^{r_2} = 0$ and $\zeta_2^{r_2} = \deg(E_2)F \cdot \zeta_2^{r_2 - 1}$ proves our case.
\end{proof}
Now we are ready to treat the case where both $E_1$ and $E_2$ are semistable vector bundles over $C$.
  \begin{thm}
  Let $E_1$ and $E_2$ be two semistable vector bundles over $C$ of rank $r_1$ and $r_2$ respectively with $r_1 \leq r_2$ and $X = \mathbb{P}(E_1) \times_C \mathbb{P}(E_2)$.
Then for all $k \in \{1, 2, ..., r_1 + r_2 - 1 \}$

$$
\overline{\Eff}^k(X) = 
\begin{cases}
\Bigg\langle \Big\{ (\zeta_1 - \mu_1F)^i (\zeta_2 - \mu_2F)^{k - i} \Big\}_{i = 0}^k, \Big\{ F \cdot \zeta_1^j \cdot\zeta_2^{k - j - 1} \Big\}_{j = 0}^{k - 1} \Bigg\rangle & if \quad k< r_1 \\ \\
\Bigg\langle \Big\{ (\zeta_1 - \mu_1F)^i (\zeta_2 - \mu_2F)^{k - i} \Big\}_{i = 0}^{r_1 - 1}, \Big\{ F \cdot \zeta_1^j \cdot \zeta_2^{k - j - 1} \Big\}_{j = 0}^ {r_1 - 1} \Bigg\rangle & if \quad r_1 \leq k < r_2 \\ \\

\Bigg\langle \Big\{ (\zeta_1 - \mu_1F)^i (\zeta_2 - \mu_2F)^{k - i} \Big\}_{i = t +1}^ {r_1 - 1}, \Big\{ F \cdot \zeta_1^j \cdot \zeta_2^{k - j -1} \Big\}_{j = t}^ {r_1 - 1} \Bigg\rangle & if \quad k = r_2 + t, \quad t = 0,..., r_1-1 .
\end{cases}
$$
where $\mu_1 = \mu(E_1)$ and $\mu_2 = \mu(E_2)$.
  \end{thm}
  \vspace{4mm}
  \begin{proof}
  Firstly, $(\zeta_1 - \mu_1F)^i \cdot(\zeta_2 - \mu_2F)^{k - i}$ and $ F\cdot \zeta_1^i \cdot \zeta_2^{k - j - 1} [ = F \cdot (\zeta_1 - \mu_1F)^i \cdot(\zeta_2 - \mu_2F)^{k - j -1} ]$ are intersections of nef divisors. So, they are pseudo-effective for all $i \in\{0, 1, 2, ..., k \}$.
 conversely, when $ k < r_1$ \, notice that we can write any element $C$ of $\overline{\Eff}^k(X)$ as \begin{align*}
  C = \sum_{i = 0}^k a_i(\zeta_1 - \mu_1F)^i \cdot (\zeta_2 - \mu_2F)^{k - i} + \sum_{j = 0}^k b_j F\cdot\zeta_1^j \cdot\zeta_2^{k - j -1}
\end{align*}  
where $a_i, b_i \in \mathbb{R}$.

 For a fixed $i_1$ intersect $C$ with $D_{i_1} :=F \cdot(\zeta_1 - \mu_1F)^{r_1 - i_1 - 1} \cdot(\zeta_2 - \mu_2F)^{r_2  - k + i_1 -1}$ and for a fixed $j_1$ intersect $C$ with $D_{j_1}:= (\zeta_1 - \mu_1F)^{r_1 - j_1 - 1} \cdot(\zeta_2 - \mu_2F)^{r_2 +j_1 - k}$. These intersections lead us to 
\begin{align*}
C \cdot D_{i_1} = a_{i_1}\quad and \quad C\cdot D_{j_1} = b_{j_1}
\end{align*}
Since $C \in \overline{\Eff}^k(X)$ and $D_{i_1}, D_{j_1}$ are intersection of nef divisors, $a_{i_1}$ and $b_{j_1}$ are non-negetive. Now running $i_1$ and$j_1$ through $\{ 0, 1, 2, ..., k \}$ we get all the $a_i$'s and $b_i$'s are non-negetive and that proves our result for $k < r_1$. The cases where $r_1 \leq k < r_2$ and $k \geq r_2$ can be proved very similarly after the intersection products involving $\zeta_1$ and $\zeta_2$ in page $2$ are taken into count.
  \end{proof}
  \bigskip
 Next we study the more interesting case where $E_1$ and $E_2$ are two unstable vector bundles of rank $r_1$ and $r_2$  and degree $d_1$ and $d_2$ respectively over a smooth curve $C$.

\smallskip
Let $E_1$ be the unique Harder-Narasimhan filtration 
\begin{align*}
E_1 = E_{10} \supset E_{11} \supset ... \supset E_{1l_1} = 0
\end{align*}
with $Q_{1i} := E_{1(i-1)}/ E_{1i}$ being semistable for all $i \in [1,l_1-1]$. Denote $ n_{1i} = \rank(Q_{1i}), \\
d_{1i} = \deg(Q_{1i})$ and $\mu_{1i} = \mu(Q_{1i}) := \frac{d_{1i}}{n_{1i}}$ for all $i$.

Similarly,  $E_2$ also admits the unique Harder-Narasimhan filtration 
\begin{align*}
E_2 = E_{20} \supset E_{21} \supset ... \supset E_{2l_2} = 0
\end{align*}
with $ Q_{2i} := E_{2(i-1)} / E_{2i}$ being semistable for $i \in [1,l_2-1]$. Denote $n_{2i} = \rank(Q_{2i}), \\ 
d_{2i} = \deg(Q_{2i})$ and $\mu_{2i} = \mu(Q_{2i}) := \frac{d_{2i}}{n_{2i}}$ for all $i$.

Consider the natural inclusion  $ \overline{i} = i_1 \times i_2 : \mathbb{P}(Q_{11}) \times_C \mathbb{P}(Q_{21}) \longrightarrow \mathbb{P}(E_1) \times_C \mathbb{P}(E_2)$, which is induced by natural inclusions $i_1 : \mathbb{P}(Q_{11}) \longrightarrow \mathbb{P}(E_1)$ and $i_2 : \mathbb{P}(Q_{21}) \longrightarrow \mathbb{P}(E_2)$. In the next theorem we will see that the cycles of $ \mathbb{P}(E_1) \times_C \mathbb{P}(E_2)$ of dimension at most $n_{11} + n_{21} - 1$ can be tied down to cycles of $\mathbb{P}(Q_{11}) \times_C \mathbb{P}(Q_{21})$ via $\overline{i}$.

  \begin{thm}
  Let $E_1$ and $E_2$ be two unstable bundle of rank $r_1$ and $r_2$ and degree $d_1$ and $d_2$ respectively over a smooth curve $C$ and $r_1 \leq r_2$ without loss of generality and $X = \mathbb{P}(E_1) \times_C \mathbb{P}(E_2)$.

\vspace{2.5mm}
Then for all $k \in \{1, 2, ..., \mathbf{n} \} \, \, (\mathbf{n} := n_{11} + n_{21} - 1)$

\smallskip
$Case(1)$: \quad $n_{11} \leq n_{21}$

$$
\overline{\Eff}_k(X) =
\begin{cases}
\Bigg\langle\Big\{ [\mathbb{P}(Q_{11}) \times \mathbb{P}(Q_{21})] (\zeta_1 - \mu_{11}F)^i (\zeta_2 - \mu_{21}F)^{\mathbf{n} - k - i} \Big\}_{i = t + 1}^{n_{11} - 1}, \Big\{ F \cdot \zeta_1^{r_1 - n_{11} +j} \cdot \zeta_2^{r_2 + n_{11} - k - j - 2} \Big\}_{j = t}^{n_{11} - 1} \Bigg\rangle & \\ \qquad \qquad if \quad k < n_{11} \quad and \quad t = 0, 1, 2, ..., n_{11} - 2 \\ \\
\Bigg\langle\Big\{ [\mathbb{P}(Q_{11}) \times \mathbb{P}(Q_{21})] (\zeta_1 - \mu_{11}F)^i (\zeta_2 - \mu_{21}F)^{\mathbf{n} - k - i}\Big\}_{i = 0}^{n_{11} - 1}, \Big\{ F\cdot \zeta_1^{r_1 - n_{11} + j} \cdot \zeta_2^{r_2 + n_{11} - k - j - 2} \Big\}_{j = 0}^{n_{11} - 1} \Bigg\rangle & \\ \qquad \qquad if \quad n_{11} \leq k < n_{21}. \\ \\

\Bigg\langle\Big\{ [\mathbb{P}(Q_{11}) \times \mathbb{P}(Q_{21})] (\zeta_1 - \mu_{11}F)^i (\zeta_2 - \mu_{21}F)^{ \mathbf{n} - k - i}\Big\}_{i = 0}^{\mathbf{n} - k}, \Big\{ F \cdot \zeta_1^{r_1 - n_{11} + j} \cdot \zeta_2^{r_2 + n_{11} - k - j - 2} \Big\}_{j = 0}^{ \mathbf{n} - k} \Bigg\rangle & \\ \qquad \qquad if \quad k \geq n_{21}.
\end{cases}
$$
\bigskip
$Case(2)$: \quad $n_{21} \leq n_{11}$

$$
\overline{\Eff}_k(X) =
\begin{cases}
\Bigg\langle\Big\{ [\mathbb{P}(Q_{11}) \times \mathbb{P}(Q_{21})] (\zeta_2 - \mu_{21}F)^i (\zeta_1 - \mu_{11}F)^{\mathbf{n} - k - i} \Big\}_{i = t + 1}^{n_{21} - 1}, \Big\{ F \cdot \zeta_2^{r_2 - n_{21} +j} \cdot \zeta_1^{r_1 + n_{21} - k - j - 2} \Big\}_{j = t}^{n_{21} - 1} \Bigg\rangle & \\ \qquad \qquad if \quad k < n_{21} \quad and \quad t = 0, 1, 2, ..., n_{21} - 2 \\ \\
\Bigg\langle\Big\{ [\mathbb{P}(Q_{11}) \times \mathbb{P}(Q_{21})] (\zeta_2 - \mu_{21}F)^i (\zeta_1 - \mu_{11}F)^{\mathbf{n} - k - i}\Big\}_{i = 0}^{n_{21} - 1}, \Big\{ F\cdot \zeta_2^{r_2 - n_{21} + j} \cdot \zeta_1^{r_1 + n_{21} - k - j - 2} \Big\}_{j = 0}^{n_{21} - 1} \Bigg\rangle & \\ \qquad \qquad if \quad n_{21} \leq k < n_{11}. \\ \\
\Bigg\langle\Big\{ [\mathbb{P}(Q_{11}) \times \mathbb{P}(Q_{21})] (\zeta_2 - \mu_{21}F)^i (\zeta_1 - \mu_{11}F)^{ \mathbf{n} - k - i}\Big\}_{i = 0}^{\mathbf{n} - k}, \Big\{ F \cdot \zeta_2^{r_2 - n_{21} + j} \cdot \zeta_1^{r_1 + n_{21} - k - j - 2} \Big\}_{j = 0}^{ \mathbf{n} - k} \Bigg\rangle & \\ \qquad \qquad if \quad k \geq n_{11}.
\end{cases}
$$
Thus in both cases $ \overline{i}_\ast$ induces an isomorphism between $ \overline{\Eff}_k([\mathbb{P}(Q_{11}) \times_C \mathbb{P}(Q_{21})])$ and $\overline{\Eff}_k(X)$ for $ k \leq \mathbf{n}$.
  \end{thm}
  \bigskip
\begin{proof}
to begin with consider $Case(1)$ and then take $ k \geq n_{21}$. Since $ (\zeta_1 - \mu_{11}F)$ and $ (\zeta_2 - \mu_{21}F)$ are nef
\begin{align*}
\phi_i := [\mathbb{P}(Q_{11}) \times \mathbb{P}(Q_{21})] (\zeta_1 - \mu_{11}F)^i (\zeta_2 - \mu_{21}F)^{ \mathbf{n} - k - i} \in \overline{\Eff}_k(X).
\end{align*}
for all $i \in \{ 0, 1, 2, ..., \mathbf{n} -k \}$.

Now The result in [Example 3.2.17, \cite{Ful}] adjusted to bundles of quotients over curves shows that
\begin{align*}
[\mathbb{P}(Q_{11})] = \eta_1^{r_1 - n_{11}} + (d_{11} - d_1)\eta_1^{r_1 - n_{11} - 1}f_1
\end{align*}
and 
\begin{align*}
[\mathbb{P}(Q_{21})] = \eta_2^{r_2 - n_{21}} + (d_{21} - d_2)\eta_2^{r_2 - n_{21} - 1}f_2
\end{align*}
Also, $p_1^\ast[\mathbb{P}(Q_{11})] \cdot p_2^\ast[\mathbb{P}(Q_{21})] = [\mathbb{P}(Q_{11}) \times_C \mathbb{P}(Q_{21})]$ . With little calculations it can be shown that

$\phi_i \cdot (\zeta_1 - \mu_{11}F)^{n_{11} - i} \cdot (\zeta_2 - \mu_{21}F)^{k + i + 1 - n_{11}}$ \\
$ =  (\zeta_1^{r_1 - n_{11}} + (d_{11} - d_1)F \cdot \zeta_1^{r_1 - n_{11} - 1})(\zeta_2^{r_2 - n_{21}} + (d_{21} - d_2)F \cdot \zeta_2^{r_2 - n_{21} - 1})(\zeta_1 - \mu_{11}F)^{n_11 - i} \cdot (\zeta_2 - \mu_{21}F)^{k + i + 1 - n_{11}}$  

$= (\zeta_1^{r_1} - d_1F \cdot \zeta_1^{r_1 - 1})(\zeta_2^{r_2} - d_2F\cdot \zeta_2^{r_2 - 1})$
 $= 0$.
 
 So, $\phi_i$ 's are in the boundary of $\overline{\Eff}_k(X)$ for all $ i \in \{0, 1, ..., \mathbf{n} - k \}$. The fact that $F \cdot \zeta_1^{r_1 - n_{11} + j} \cdot \zeta_2^{r_2 + n_{11} - k - j - 2}$ 's are in the boundary of $\overline{\Eff}_k(X)$ for all $ i \in \{0, 1, ..., \mathbf{n} - k \}$ can be deduced from the proof of Theorem 2.2. The other cases can be proved similarly.
 
 The proof of $Case(2)$ is similar to the proof of $Case(1)$.
 
 Now, to show the isomorphism between pseudo-effective cones induced by $\overline{i}_\ast$ observe that $Q_{11}$ and $Q_{21}$ are semi-stable bundles over $C$. So, Theorem 2.2 gives the expressions for  $\overline{\Eff}_k([\mathbb{P}(Q_{11}) \times_C \mathbb{P}(Q_{21})])$. Let $\zeta_{11} = \mathcal{O}_{\mathbb{P}(\tilde{\pi}_2^\ast(Q_{11}))}(1) =\tilde{p}_1^\ast(\mathcal{O}_{\mathbb{P}(Q_{11})}(1))$  and $\zeta_{21} = \mathcal{O}_{\mathbb{P}(\tilde{\pi}_1^\ast(Q_{21})}(1) = \tilde{p}_2^\ast(\mathcal{O}_{\mathbb{P}(Q_{21})}(1))$, where $\tilde{\pi_2} = \pi_2|_{\mathbb{P}(Q_{21})}$, $ \tilde{\pi_1} = \pi_1|_{\mathbb{P}(Q_{11})}$ and $\tilde{p}_1 : \mathbb{P}(Q_{11}) \times_C \mathbb{P}(Q_{21}) \longrightarrow \mathbb{P}(Q_{11})$, $\tilde{p}_2 : \mathbb{P}(Q_{11}) \times_C\mathbb{P}(Q_{21}) \longrightarrow \mathbb{P}(Q_{21})$ are the projection maps. Also notice that $ \overline{i}^\ast \zeta_1 = \zeta_{11}$ and $\overline{i}^\ast \zeta_1 = \zeta_{21}$.
 
 Using the above relations and projection formula the isomorphism between $ \overline{\Eff}_k([\mathbb{P}(Q_{11}) \times_C \mathbb{P}(Q_{21})])$ and $\overline{\Eff}_k(X)$ for $ k \leq \mathbf{n}$ can be proved easily.

\end{proof}

Next we want to show that higher dimension pseudo effective cycles on $X$ can be related to the pseudo effective cycles on $ \mathbb{P}(E_{11}) \times_C \mathbb{P}(E_{21})$. More precisely there is a isomorphism between $\overline{\Eff}^k(X)$ and $\overline{\Eff}^k([\mathbb{P}(E_{11}) \times_C \mathbb{P}(E_{21})])$ for $ k < r_1 + r_2 - 1 - \mathbf{n}$. Useing the coning construction as in [ful] we show this in two steps, first we establish an isomorphism between $\overline{\Eff}^k([\mathbb{P}(E_1) \times_C \mathbb{P}(E_2)])$ and $\overline{\Eff}^k([\mathbb{P}(E_{11}) \times_C \mathbb{P}(E_2)])$ and then an isomorphism between $\overline{\Eff}^k([\mathbb{P}(E_{11}) \times_C \mathbb{P}(E_2)])$ and $\overline{\Eff}^k([\mathbb{P}(E_{11}) \times_C \mathbb{P}(E_{21})])$ in similar fashion. But before proceeding any further we need to explore some more facts.

Let $E$ be an unstable vector bundle over a non-singular projective variety $V$. There is a unique filtration 
\begin{align*}
E = E^0 \supset E^1 \supset E^2 \supset ... \supset E^l = 0
\end{align*}
which is called the harder-Narasimhan filtration of $E$ with $Q^i := E^{i - 1}/E^i$ being semistable for $ i \in [1, l - 1]$. Now the following short-exact sequence
\begin{align*}
0 \longrightarrow E^1 \longrightarrow E \longrightarrow Q^1 \longrightarrow 0
\end{align*}
induced by the harder-narasimhan filtration of $E$ gives us the natural inclusion $j : \mathbb{P}(Q^1) \hookrightarrow \mathbb{P}(E)$. Considering $\mathbb{P}(Q^1)$ as a subscheme of $\mathbb{P}(E)$ we obtain the commutative diagram below by blowing up $\mathbb{P}(Q^1)$.

\begin{center}
\begin{equation}
\begin{tikzcd}
\tilde{Y} = Bl_{\mathbb{P}(Q^1)}{\mathbb{P}(E)} \arrow[r, "\Phi"] \arrow[d, "\Psi"] & \mathbb{P}(E^1) = Z \arrow [d, "q"] \\
Y = \mathbb{P}(E) \arrow [r, "p"] 
& V
\end{tikzcd}
\end{equation}
\end{center}
where $\Psi$ is blow-down map.
\begin{thm}
With the above notation,there exists a locally free sheaf $G$ on $Z$ such that $\tilde{Y} \simeq \mathbb{P}_Z(G)$ and $\nu : \mathbb{P}_Z(G) \longrightarrow Z$ it's corresponding bundle map.

In particular if we place $V = \mathbb{P}(E_2)$, $ E = \pi_2^\ast E_1$, $E^1 = \pi_2^\ast E_{11}$ and $ Q^1 = \pi_2^\ast Q_{11}$ then the above commutative diagram becomes

\begin{center}
\begin{equation}
\begin{tikzcd}
\tilde{Y'} = Bl_{\mathbb{P}(\pi_2^\ast Q_{11})}{\mathbb{P}(\pi_2^\ast E_1)} \arrow[r, "\Phi'"] \arrow[d, "\Psi'"] &
\mathbb{P}(\pi_2^\ast E_{11}) = Z' \arrow[d, "\overline{p}_2"] \\
Y' = \mathbb{P}(\pi_2^\ast E_1) \arrow[r, "p_2"]
& \mathbb{P}(E_2)
\end{tikzcd}
\end{equation}
\end{center}
where $p_2 : \mathbb{P}(\pi_2^\ast E_1) \longrightarrow \mathbb{P}(E_2)$ and $\overline{p}_2 : \mathbb{P}(\pi_2^\ast E_{11}) \longrightarrow \mathbb{P}(E_2)$ are projection maps.

and there exists a locally free sheaf $G'$ on $Z'$ such that  $\tilde{Y'} \simeq \mathbb{P}_{Z'}(G')$ and $\nu' : \mathbb{P}_{Z'}(G') \longrightarrow Z'$ it's bundle map.

 Now let $\zeta_{Z'} = \mathcal{O}_{Z'}(1)$, $\gamma = \mathcal{O}_{\mathbb{P}_{Z'}(G')}(1)$, $F$ the numerical equivalence class of a fibre of  $\pi_2 \circ p_2$, $F_1$ the numerical equivalence class of a fibre of  $\pi_2 \circ \overline{p}_2$, $\tilde{E}$ the class of the exceptional divisor of $\Psi'$ and $\zeta_1 = p_1^\ast(\eta_1) = \mathcal{O}_{\mathbb{P}(\pi_2^\ast E_1)}(1)$. Then we have the following relations:
 \vspace{4mm}
 \begin{align}
 \gamma = (\Psi')^\ast \, \zeta_1, \quad (\Phi')^\ast \, \zeta_{Z'} = (\Psi')^\ast \, \zeta_1 - \tilde{E}, \quad (\Phi')^\ast F_1 = (\Psi')^\ast F 
 \end{align}
 \begin{align}
 \tilde{E} \cdot (\Psi')^\ast \, (\zeta_1 - \mu_{11}F)^{n_{11}} = 0
 \end{align}
  
 Additionaly, if we also denote the support of the exceptional divisor of $\tilde{Y'}$ by $\tilde{E}$ , then $\tilde{E} \cdot N(\tilde{Y'}) = (j_{\tilde{E}})_\ast N(\tilde{E})$, where  $j_{\tilde{E}}: \tilde{E} \longrightarrow \tilde{Y'}$ is the canonical inclusion.
\end{thm}
\vspace{4mm}
\begin{proof}
With the above hypothesis the following commutative diagram is formed:

\begin{center}
\begin{tikzcd}%
0 \arrow[r] & q^\ast E^1  \arrow[r] \arrow[d, two heads] & q^\ast E  \arrow[r] \arrow[d] & q^\ast Q^1  \arrow[r] \ar[equal] {d} & 0 \\
0 \arrow[r] & \mathcal{O}_{\mathbb{P}(E^1)}(1) \arrow[r] & G \arrow[r] & q^\ast Q^1 \arrow[r] & 0
\end{tikzcd}%
\end{center}
where $G$ is the push-out of morphisms $ q^\ast E^1 \longrightarrow q^\ast E$ and $ q^\ast E^1 \longrightarrow \mathcal{O}_{\mathbb{P}(E^1)}(1)$ and the first vertical map is the natural surjection. Now let $W = \mathbb{P}_Z(G)$ and  $\nu : W \longrightarrow Z$ be it's bundle map. So there is a cannonical surjection $\nu^\ast G \longrightarrow \mathcal{O}_{\mathbb{P}_Z(G)}(1)$. Also note that $q^\ast E \longrightarrow G$ is surjective by snake lemma. Combining these two we obtain a surjective morphism $\nu^\ast q^\ast E \longrightarrow \mathcal{O}_{\mathbb{P}_Z(G)}(1)$ which determines $\omega : W \longrightarrow Y$. We claim that we can identify $(\tilde{Y}, \Phi, \Psi)$ and $(W, \nu, \omega)$. Now Consider the following commutative diagram:

\begin{center}
\begin{equation}
\begin{tikzcd}
W = \mathbb{P}_Z(G)
\arrow[drr, bend left, "\nu"]
\arrow[ddr, bend right, "\omega"]
\arrow[dr, "\mathbf{i}"] & & \\
& Y \times_V Z = \mathbb{P}_Z(q^\ast E) \arrow[r, "pr_2"] \arrow[d, "pr_1"]
& \mathbb{P}(E^1) = Z \arrow[d, "q"] \\
& Y = \mathbb{P}(E) \arrow[r, "p"]
& V
\end{tikzcd}
\end{equation}
\end{center}
\vspace{12mm}
 where $\mathbf{i}$ is induced by the universal property of the fiber product. Since $\mathbf{i}$ can also be obtained from the surjective morphism $q^\ast E \longrightarrow G$ it is a closed immersion. Let $\mathcal{T}$ be the $\mathcal{O}_Y$ algebra $\mathcal{O}_Y \oplus \mathcal{I} \oplus \mathcal{I}^2 \oplus ...$, where $\mathcal{I}$ is the ideal sheaf of $\mathbb{P}(Q^1)$ in $Y$. We have an induced map of $\mathcal{O}_Y$- algebras $Sym(p^\ast E^1) \longrightarrow \mathcal{T} \ast \mathcal{O}_Y(1)$ which is onto because the image of the composition $ p^\ast E^1 \longrightarrow p^\ast E \longrightarrow \mathcal{O}_Y(1)$ is  $ \mathcal{T} \otimes \mathcal{O}_Y(1)$. This induces a closed immersion
 \begin{center}
  $\mathbf{i}' : \tilde{Y} = Proj(\mathcal{T} \ast \mathcal{O}_Y(1)) \longrightarrow Proj(Sym(p^\ast E^1) = Y \times_V Z$.
  \end{center}
 $\mathbf{i}'$ fits to a similar commutative diagram as $(5)$ and as a result $\Phi$ and $\Psi$ factor through $pr_2$ and $pr_1$.
 Both $W$ and $\tilde{Y}$ lie inside $Y \times_V Z$ and $\omega$ and $\Psi$ factor through $pr_1$ and $\nu$ and $\Phi$ factor through $pr_2$. So to prove the identification between $(\tilde{Y}, \Phi, \Psi )$ and $(W, \nu, \omega)$ , it is enough to show that $ \tilde{Y} \cong W$. This can be checked locally. So, after choosing a suitable open cover for $V$ it is enough to prove  $\tilde{Y} \cong W$ restricted to each of these open sets. Also we know that $p^{-1}(U) \cong \mathbb{P}_U ^{rk(E) - 1}$ when $E_{|U}$ is trivial and $\mathbb{P}_U^n = \mathbb{P}_{\mathbb{C}}^n \times U$. Now the the isomorphism follows from [proposition 9.11, \cite{EH}] after adjusting the the definition of projectivization in terms of \cite{H}.
 
 We now turn our attention to the diagram $(3)$. observe that if we fix the notations $W' = \mathbb{P}_Z'(G')$ with $\omega' : W' \longrightarrow Y'$ as discussed above then we have an identification between $(\tilde{Y}', \Phi', \Psi')$ and $(W', \nu', \omega')$.
 
 $\omega' : W' \longrightarrow Y'$  comes with $(\omega')^\ast \mathcal{O}_
 {Y'}(1) = \mathcal{O}_{\mathbb{P}_{Z'}(G')}(1)$. So, $\gamma = (\Psi')^\ast \, \zeta_1$ is achieved. $(\Phi')^\ast F_1 = (\Psi')^\ast F$  follows from the commutativity of the diagram $(3)$. 
 
 The closed immersion $\mathbf{i}'$ induces a relation between the $\mathcal{O}(1)$ sheaves of $Y \times_V Z$ and $\tilde{Y}$. For $Y \times_V Z$ the $\mathcal{O}(1)$ sheaf is $pr_2^ \ast \mathcal{O}_{Z}(1)$ and for $ Proj(\mathcal{T} \ast \mathcal{O}_Y(1)$ the $\mathcal{O}(1)$ sheaf is $\mathcal{O}_{\tilde{Y}}( - \tilde{E}) \otimes (\Psi)^ \ast \mathcal{O}_Y(1)$. Since $\Phi$ factors through $pr_2$, $(\Phi)^ \ast \mathcal{O}_Z(1) = \mathcal{O}_{\tilde{Y}}( - \tilde{E}) \otimes (\Psi)^ \ast \mathcal{O}_Y(1)$. In the particular case (see diagram $(3)$) $(\Phi')^ \ast \mathcal{O}_{Z'}(1) = \mathcal{O}_{\tilde{Y}'}( - \tilde{E}) \otimes (\Psi')^ \ast \mathcal{O}_{Y'}(1)$ i. e. $(\Phi')^\ast \, \zeta_{Z'} = (\Psi')^\ast \, \zeta_1 - \tilde{E} $.
 
 Next consider the short exact sequence: 
 \begin{center}$ 0 \longrightarrow \mathcal{O}_{Z'}(1) \longrightarrow G' \longrightarrow \overline{p}_2^\ast \pi_2^ \ast Q_{11} \longrightarrow 0$
 \end{center}
  We wish to calculate below the total chern class of $G'$ through the chern class relation obtained from the above short exact sequence.
  \begin{center}
  $c(G') = c(\mathcal{O}_{Z'}(1)) \cdot c(\overline{p}_2^\ast \pi_2^ \ast Q_{11}) = (1 + \zeta_{Z'}) \cdot \overline{p}_2^\ast \pi_2^\ast(1 + d_{11}[pt]) = (1 + \xi_{Z'})(1 + d_{11}F_1)$
  \end{center}
 From the grothendieck relation for $G'$ we have
 
 \bigskip
 
 $\gamma^{n_{11} + 1} - {\Phi'} ^ \ast(\zeta_{Z'} + d_{11}F_1) \cdot \gamma^{n_{11}} + {\Phi'} ^ \ast (d_{11}F_1 \cdot \zeta_{Z'}) \cdot \gamma^ {n_{11} - 1} = 0$ \\
 \vspace{1mm}
 $\Rightarrow  \gamma^{n_{11} + 1} - ({\Psi'} ^ \ast \zeta_1 - \tilde{E}) + d_{11}{\Psi'}^ \ast F) \cdot \gamma^{n_{11}} + d_{11}({\Psi'} ^ \ast \zeta_1 - \tilde{E}) \cdot {\Psi'}^ \ast F) \cdot \gamma^{n_{11} - 1} = 0$ \\
 \vspace{1mm}
 $\Rightarrow \tilde{E} \cdot \gamma^{n_{11}} - d_{11}\tilde{E} \cdot {\Psi'} ^ \ast F \cdot \gamma^{n_{11} - 1} = 0$\\
 \vspace{1mm}
 $\Rightarrow \tilde{E} \cdot {\Psi'}^ \ast (\zeta_1 - \mu_{11}F)^{n_{11}} = 0$
 
 For the last part note that $\tilde{E} = \mathbb{P}(\pi_2^\ast Q_{11}) \times_{\mathbb{P}(E_2)} Z'$. Also  $N(\tilde{Y}')$ and $N(\tilde{E})$ are free $N(Z')$-module. Using these informations and projection formula, the identity $\tilde{E} \cdot N(\tilde{Y'}) = (j_{\tilde{E}})_\ast N(\tilde{E})$ is obtained easily.

\end{proof}
Now we are in a position to prove the next theorem.
\begin{thm}
$\overline{\Eff}^k(X) \cong \overline{\Eff}^k(Y') \cong \overline{\Eff}^k(Z')$ and $\overline{\Eff}^k(Z') \cong \overline{\Eff}^k(Z'')$. So, $\overline{\Eff}^k(X) \cong \overline{\Eff}^k(Z'')$ for $k < r_1 + r_2 - 1 - \mathbf{n}$

where $Z'= \mathbb{P}(E_{11}) \times_C \mathbb{P}(E_2)$ and $ Z'' = \mathbb{P}(E_{11}) \times_C \mathbb{P}(E_{21})$
\end{thm}
\vspace{4mm}
\begin{proof}
Since $Y' = \mathbb{P}(\pi_2^ \ast E_1) \cong \mathbb{P}(E_1) \times_C \mathbb{P}(E_2) = X$, \, $\overline{\Eff}^k(X) \cong \overline{\Eff}^k(Y')$ is followed at once. To prove that $\overline{\Eff}^k(X) \cong \overline{\Eff}^k(Z')$ we first define the the map:
$\theta_k: N^k(X) \longrightarrow N^k(Z')$
by
\begin{align*}
\zeta_1^ i \cdot \zeta_2^ {k - i} \mapsto \bar{\zeta_1}^  i \cdot \bar{\zeta_2}^{k - i}, \quad F\ \cdot \zeta_1^j \cdot \zeta_2^{k - j - 1} \mapsto F_1 \cdot \bar{\zeta}_1^ j \cdot \bar{\zeta}_2^ {k - j - 1}
\end{align*}
\smallskip
where $\bar{\zeta_1} = \overline{p}_1 ^\ast(\mathcal{O}_{\mathbb{P}(E_{11})}(1))$ and $\bar{\zeta_2} = \overline{p}_2 ^\ast(\mathcal{O}_{\mathbb{P}(E_2)}(1))$. $ \overline{p}_1 : \mathbb{P}(E_{11}) \times_C \mathbb{P}(E_2) \longrightarrow \mathbb{P}(E_{11})$ and $ \overline{p}_2 : \mathbb{P}(E_{11}) \times_C \mathbb{P}(E_2) \longrightarrow \mathbb{P}(E_2)$ are respective projection maps.

It is evident that the above map is in isomorphism of abstract groups. We claim that this induces an isomorphism between $ \overline{\Eff}^k(X)$ and $\overline{\Eff}(Z')$. First we construct an inverse for $\theta_k$. Define $\Omega_k : N^k(Z') \longrightarrow N^k(X)$ by
\begin{center}
$\Omega_k (l) = {\Psi'}_\ast {\Phi'}^\ast (l)$
\end{center}
$\Omega_k$ is well defined since $\Phi'$ is flat and $\Psi'$ is birational. $\Omega_k$ is also pseudo-effective. Now we need to show that $\Omega_k$ is the inverse of $\theta_k$.

\begin{align*}
\Omega_k(\bar{\zeta_1}^  i \cdot \bar{\zeta_2}^{k - i}) &  = {\Psi'}_\ast (({\Phi'}^ \ast \bar{\zeta_1})^i \cdot ({\Phi'}^\ast \bar{\zeta_2})^{k - i}) \\
& = {\Psi'}_\ast (({\Phi'}^\ast \zeta_{Z'})^i \cdot ({\Phi'}^\ast \bar{\zeta_2})^{k - i}) \\
& = {\Psi'}_\ast (({\Psi'}^\ast \zeta_1 - \tilde{E})^i \cdot ({\Psi'}^\ast \zeta_2)^{k - i}) \\
& = {\Psi'}_\ast((\sum_{0 \leq c \leq i} (-1)^ i \tilde{E}^c ({\Psi'}^\ast \zeta_1)^{i - c}) \cdot ({\Psi'}^\ast \zeta_2)^{k - i}) \\
\end{align*}
\smallskip
Similarly,
\begin{align*}
\Omega_k(F_1 \cdot \bar{\zeta_1}^j \cdot \bar{\zeta_2}^{k - j -1}) = {\Psi'}_\ast (( \sum_{0 \leq d \leq j} (-1)^j \tilde{E}^d ({\Psi'}^\ast \zeta_1)^{j - d}) \cdot ({\Psi'}^\ast \zeta_2)^{k - j - 1})
\end{align*}

So,

$\Omega_k\Big(\sum_i a_i\, \bar{\zeta_1}^  i \cdot \bar{\zeta_2}^{k - i} + \sum_j b_j \, F_1 \cdot \bar{\zeta_1}^j \cdot \bar{\zeta_2}^{k - j -1} \Big)$
\begin{align*}
 = \Big(\sum_i a_i\,{\zeta_1}^  i \cdot{\zeta_2}^{k - i} + \sum_j b_j \, F \cdot{\zeta_1}^j \cdot{\zeta_2}^{k - j -1} \Big) + {\Psi'}_\ast \Big( \sum_i \sum_{1 \leq c \leq i} \tilde{E}^c {\Psi'}^\ast(\alpha_{i, c}) + \sum_j \sum_{1 \leq d \leq j} \tilde{E}^d {\Psi'}^\ast (\beta_{j, d})\Big)
\end{align*}
for some cycles $\alpha_{i, c}, \beta_{j, d} \in N(X)$.
But, ${\Psi'}^\ast(\tilde{E}^t) = 0$ for all $1 \leq t \leq i \leq r_1 + r_2  - 1 - \mathbf{n}$ for dimensional reasons. Hence, the second part in the right hand side of the above equation vanishes and we make the conclusion that $\Omega_k = \theta_k ^{-1}$.

Next we seek an inverse of $\Omega_k$ which is pseudo-effective and meet our demand of being equal to $\theta_k$. Define $\eta_k : N^k(X) \longrightarrow N^k(Z')$ by
\begin{align*}
\eta_k(s) = {\Phi'}_\ast(\delta \cdot {\Psi'}^\ast s)
\end{align*}
where $ \delta = {\Psi'}^\ast (\xi_2 - \mu_{11}F)^{n_{11}}$.

By the relations $(5)$ and $(6)$  , ${\Psi'}^ \ast( (\zeta_1^i \cdot \zeta_2^{k - i})$ is ${\Phi'}^\ast(\bar{\zeta_1}^i \cdot \bar{\zeta_2}^{k - i})$ modulo $\tilde{E}$ and $\delta \cdot \tilde{E} = 0$. Also ${\phi'}_\ast \delta = [Z']$ which is derived from the fact that ${\Phi'}_\ast \gamma^{n_{11}} = [Z']$ and the same relations $(5)$ and $(6)$.  Therefore 
\smallskip
\begin{align*}
\eta_k(\zeta_1^i \cdot \zeta_2^{k - i}) =  {\Phi'}_\ast (\delta \cdot {\Phi'}^\ast(\bar{\zeta_1}^i \cdot \bar{\zeta_2}^{k - i})) = (\bar{\zeta_1}^i \cdot \bar{\zeta_2}^{k - i}) \cdot [Z'] = \bar{\zeta_1}^i \cdot \bar{\zeta_2}^{k - i}
\end{align*}
In a similar way, ${\Psi'}^ \ast(F \cdot \zeta_1^ j \cdot \zeta_2^{k - j - 1})$ is ${\Phi'}^\ast (F_1 \cdot \bar{\zeta_1}^j \cdot \bar{\zeta_2}^{k - j -1})$ modulo $\tilde{\mathbf{E}}$ and as a result of this
\begin{align*}
\eta_k(F \cdot \zeta_1^ j \cdot \zeta_2^{k - j - 1}) = F_1 \cdot \bar{\zeta_1}^j \cdot \bar{\zeta_2}^{k - j -1}
\end{align*}

So, $\eta_k = \theta_k$.

Next we need to show that $\eta_k$ is a pseudo- effective map. Notice that ${\Psi'}^\ast s = \bar{s} + \mathbf{j}_\ast s'$ for any effective cycle $s$ on $X$, where $\bar{s}$ is the strict transform under $\Psi'$ and hence effective. Now  $\delta$ is intersection of nef classes. So, $\delta \cdot \bar{s}$ is pseudo-effective. Also $\delta \cdot \mathbf{j}_\ast s' = 0$ from theorem 2.4 and ${\Phi'}_\ast$ is pseudo-effective. Therefore $\eta_k$ is pseudo-effective and first part of the theorem is proved. We will sketch the prove for the second part i.e. $\overline{\Eff}^k(Z') \cong \overline{\Eff}^k(Z'')$ which is similar to the proof of the first part. Consider the following diagram:

\begin{center}
\begin{equation}
\begin{tikzcd}
 Z'' = \mathbb{P}(E_{11}) \times_C \mathbb{P}(E_{21}) \arrow[r, "\hat{p}_2"] \arrow[d, "\hat{p}_1"]
& \mathbb{P}(E_{21})\arrow[d,"\hat{\pi}_2"]\\
\mathbb{P}(E_{11}) \arrow[r, "\hat{\pi}_1" ]
& C
\end{tikzcd}
\end{equation}
\end{center}
Define $\hat{\theta}_k : N^k(Z') \longrightarrow N^k(Z'')$ by 
\begin{align*}
\bar{\zeta_1}^ i \cdot \bar{\zeta_2}^  {k - i} \mapsto \hat{\zeta_1}^ i \cdot \hat{\zeta_2}^{k - i}, \quad F \cdot \bar{\zeta_1}^j \cdot \bar{\zeta_2}^ {k - j - 1} \mapsto F_2 \cdot \hat{\zeta_1}^j \cdot \hat{\zeta_2}^{k - j - 1}
\end{align*}

where $\hat{\zeta_1} = \hat{p_1}^\ast (\mathcal{O}_{\mathbb{P}(E_{11})}(1)), \hat{\zeta_2} = \hat{p_2}^\ast (\mathcal{O}_{\mathbb{P}(E_{21})}(1))$ and $F_2$ is the class of a fibre of $\hat{\pi_1} \circ \hat{p_1}$.

This is a isomorphism of abstract groups and behaves exactly the same as $\theta_k$. The methods applied to get the result for $\theta_k$ can also be applied successfully here.

\subsection*{Acknowledgement}
The author would like to thank Prof. D.S. Nagaraj, IISER Tirupati for suggestions and discussions at every stage of this work. This work is supported financially by a fellowship from IMSc,Chennai (HBNI),
DAE, Government of India.
\end{proof}
 
\end{document}